\documentclass[12pt]{article}
\def\Cliff{{\rm Cliff}}
\newtheorem{thm}{Theorem}[section]
\newtheorem{prop}[thm]{Proposition}
\newtheorem{lemma}[thm]{Lemma}
\def\proof{{\sc Proof. }}
\usepackage{amsfonts, amssymb}
\def\R{\mathbb{R}}
\def\RP{\mathbb{R}{\rm P}}
\def\C{\mathbb{C}}
\def\H{\mathbb{H}}
\def\A{\mathbb{A}}
\def\Nc{\mathcal{N}}
\def\Box{\square}
\def\<{\langle}
\def\>{\rangle}

\title{Circles and Clifford Algebras}
\author{V. Timorin
\thanks{Partially supported by CRDF RM1-2086}}
\date{}

\begin{document}
\maketitle

{\small {\bf Abstract.} Consider a smooth map from a neighborhood
of the origin in a real vector space to a neighborhood of the origin in a Euclidean
space. Suppose that this map takes all germs of lines passing
through the origin to germs of Euclidean circles, or lines, or a point. We
prove that under some simple additional assumptions this map takes
all lines passing though the origin to the same circles as a Hopf map
coming from a representation of a Clifford algebra does. We also
describe a connection between our result and the Hurwitz--Radon
theorem about sums of squares.}

\section{Introduction}

By a circle in a Euclidean space we mean an honest Euclidean
circle, or a straight line, or a point. For the sake of brevity,
by a {\em vector line} we always mean a line passing though the
origin (i.e. a 1-dimensional vector subspace).

Our problem is as follows: describe all smooth (or continuous)
maps from an open subset of an affine space $\R^m$ (or a
projective space $\RP^m$) to a Euclidean space $\R^n$ (or a sphere
$S^n$) that take germs of straight lines to germs of circles (for
short: take lines to circles). Originally, this problem came from
nomography and was stated for diffeomorphisms between open subsets
of $\R^2$ by G. S. Khovanskii \cite{GKh}. Nomography deals with
graphical representations of functions. Nomograms using a compass
are more accurate than those using a ruler. Thus the above problem
arises. The original 2-dimensional problem was solved by A.
Khovanskii in \cite{Kh}. He proved that, up to a projective
transformation in the space of pre-image and a M\"obius
transformation in the space of image, there are only 3 smooth
one-to-one maps taking lines to circles. They correspond to the 3
classical geometries --- Euclidean, spherical and hyperbolic. The
Euclidean map is just the identity. The spherical map is the
composition of a central projection to a sphere and a
stereographic projection back to the plane. The hyperbolic map is
the same with a hyperboloid instead of a sphere. This is the map
between 2 models of hyperbolic geometry: the Klein model and the
Poincar\'e unit disk model.

In a unified way all the maps from Khovanskii's theorem may be
described as follows.
To get a map from an open subset of the plane to a sphere
taking all lines to circles, one has to embed this subset
projectively to the 3-dimensional space and then project it
to a sphere from some point (possibly infinite).
If we project from inside of the sphere, then we get
the spherical case. If the center of projection lies on the
sphere, then we have the Euclidean case. Finally, projecting from
outside of the sphere, we get the hyperbolic case.

F. Izadi in his PhD thesis \cite{Iz} carried over the results of Khovanskii to
the 3-dimensional case.

But in dimension 4 all this is wrong! There are more maps taking
lines to circles. For example, identify $\R^4$ with $\C^2$ and
consider a complex projective transformation. It takes complex
lines to complex lines, and on each complex line it acts as a
M\"obius transformation. Note that every real line belongs to a
unique complex line. Thus any real line goes to a circle.

There is a large class of maps from an open subset of $\R^4$ to
$\R^4$ that take lines to circles. Let $A$ and $B$ be 2 affine
maps from $\R^4$ to $\R^4$. Identify $\R^4$ with the algebra $\H$
of quaternions. Consider the maps $x\mapsto A(x)B(x)^{-1}$ and
$x\mapsto B(x)^{-1}A(x)$. Here the multiplication and the
inversion are in the sense of quaternions. We call these maps {\em
fractional quaternionic transformations}. They take lines to
circles. Note that fractional quaternionic transformations are not
necessarily fractional linear: $A$ and $B$ are affine over reals,
not over quaternions.
In \cite{maplc4} it is proved that diffeomorphisms between open
regions in $\R^4$ taking lines to circles are either fractional
quaternionic transformations or projections from a hyperplane in
$\R^5$ to a sphere.

In \cite{Tim}, the following ``micro-local'' theorem is proved.
Consider a germ of diffeomorphism $\Phi:(\R^4,0)\to(\R^4,0)$ that
takes germs of vector lines to germs of circles. Then
the images of all vector lines are the same circles as for some
fractional quaternionic transformation preserving 0.
This describes so called
{\em rectifiable pencils of circles} (sets of circles passing
through some point and obtained from straight lines via some local
diffeomorphism). In dimensions 2 and 3, the situation is much
simpler: any sufficiently large rectifiable pencil of circles has
some other common point of intersection, i.e., it is obtained from
a pencil of lines by some M\"obius transformation \cite{Kh,Iz}.

In this paper, we study the multidimensional case. There is a
general construction of smooth maps from $\RP^m$ to $S^n$ that
take all lines to circles, the {\em Hopf maps}. The most
systematic construction of Hopf maps is based on representations
of Clifford algebras $\Cliff(r)$ generated over $\R$ by $r$
anti-commuting imaginary units. With any representation $\phi$ of
$\Cliff(r)$ in $\R^n$ we associate a map from an open subset of
$\R^{r+n+1}$ to $\R^n$. Namely, $\R^{r+1}$ is embedded into
$\Cliff(r)$ as the subspace spanned by all generators and 1. A map
$F_\phi$ from $(\R^{r+1}-0)\times\R^n$ to $\R^n$ sends a pair
$(\alpha,x)$ to $\phi(\alpha)^{-1}(x)$. This map takes lines to
circles. Projectivizing the space of preimage and compactifying
the space of image we obtain a map from $\RP^{n+r}$ to $S^n$ with
the same property. This is the Hopf map associated with $\phi$.

We will prove the following ``micro-local'' theorem. Consider a
germ of a smooth map $\Phi:(\R^r\times\R^n,0)\to(\R^n,0)$ that
takes germs of vector lines to germs of circles.
Assume that $\R^n$ is linearly embedded into $\R^{r+n}$ so that
$\Phi$ is the identity on $\R^n$. Assume also that the images of
lines not lying in $\R^n$ do not degenerate to lines. Then the image
of each vector line is the same circle as for some map $F_\phi$
or its restriction to a hyperplane.

The maps $F_\phi:\R^m\to\R^n$ give many examples of rectifiable
pencils of circles. It is enough to compose $F_\phi$ with a
generic affine embedding of $\R^n$ to $\R^m$ to get a local
diffeomorphism sending lines to circles.

Let us return to dimension 4. The algebra $\Cliff(3)$ is
isomorphic to $\H\oplus\H$. Thus it has 2 representations in
$\R^4=\H$. This gives 2 maps from $\R^8$
to $\R^4$ taking all lines to circles. The corresponding maps from
$S^7$ to $S^4$ are classical Hopf fibrations. Note that all
circles are the images of lines under these maps. All fractional
quaternionic transformations can be obtained from these maps as
described above.

The paper is organized as follows. In Section 2, we recall some basic facts
about Clifford algebras and their representations. Section 3 contains the
necessary information about the Hopf maps.
Then we formulate the micro-local theorem in Section 4. Section 5 establishes
a connection between this theorem and the Hurwitz--Radon theorem about
sums of squares. Finally, in Section 6 we conclude the proof of the micro-local
theorem and give a description of Clifford algebras representations as linear
spaces of complex multiplications.

I am grateful to A. Khovanskii for useful discussions and to the referee
for useful comments.

\section{Clifford algebras}

Let $V$ be a real vector space equipped with a symmetric bilinear form $Q$. The
Clifford algebra $\Cliff(V,Q)$ is the associative
algebra with a unit generated by all vectors from $V$ subject to relations
$vw+wv=Q(v,w)$ for all $v,w\in V$. If $V$ has dimension $r$ and $Q$
is negative definite, then the algebra $\Cliff(V,Q)$ is denoted by $\Cliff(r)$.
In other words, $\Cliff(r)$ is generated by $r$ anti-commuting imaginary units
$e_1,\dots,e_r$:
$$e_i^2=-1,\qquad e_ie_j+e_je_i=0\quad (i\ne j).$$
We will be always dealing with Clifford algebras of type $\Cliff(r)$. For
a Euclidean space $V$ with the Euclidean form $q$ we denote the algebra
$\Cliff(V,-q)$ by $\Cliff(V)$. A Clifford algebra is a direct generalization
of the skew-field of quaternions $\H=\Cliff(2)$.

The following table gives the Clifford algebras $\Cliff(r)$ for $r<8$.

\bigskip
\begin{tabular}{|c|c|c|c|c|c|c|c|c|}
\hline
$r$& 0& 1& 2& 3& 4& 5& 6& 7\\
\hline
$\Cliff(r)$& $\R$& $\C$& $\H$& $\H\oplus\H$& $\H[2]$& $\C[4]$& $\R[8]$&
$\R[8]\oplus\R[8]$\\
\hline
\end{tabular}
\bigskip

Here $\A[r]$ ($\A=\R,\C$ or $\H$) denotes the algebra of
$r\times r$-matrices over $\A$. The other Clifford algebras are computed by
the Bott periodicity
$$\Cliff(r+8)=\Cliff(r)\otimes\R[16].$$
This table was first obtained by E. Cartan \cite{Car}, see also \cite{ABS}.
All real finite-dimensional representations of Clifford algebras read from the table
since a matrix algebra $\A[r]$ has a unique irreducible representation,
namely, $\A^r$, and all representations are completely reducible.

{\bf Example 1.}
The algebra $\Cliff(3)$ is isomorphic to $\H\oplus\H$.
In this presentation,  $e_1=(i,-i)$, $e_2=(j,-j)$ and $e_3=(k,-k)$.

{\bf Example 2.}
The algebra $\Cliff(7)$ acts on $\R^8=\H^2$. Let us describe the action of the space
generated by 1 and $e_i$ ($i=1,\dots,7$). The latter is also identified with
$\H^2$. For any quaternion $a$ denote by $L_a$ the operator of left
multiplication by $a$ in $\H$ and by $R_a$ the operator of right multiplication.
A pair $(a,b)\in\H^2$ acts as
$$\left(\begin{array}{cc}L_a& -R_{\bar b}\\R_b& L_{\bar a}\end{array}\right).$$
If we interchange the left and right multiplications, then we obtain
another representation of $\Cliff(7)$ in $\R^8$.

Consider a Euclidean vector space $V$. Denote by $\hat V$ the space $V\oplus\<1\>$,
which lies in $\Cliff(V)$ and is generated by $V$ and $1$. Let us say that
a linear representation of $\Cliff(V)$ in $\R^n$ is {\em compatible with
the Euclidean metric} (or with Euclidean inner product) if
any vector from $\hat V$ acts conformally, i.e. multiplies all distances
by the same factor.
This condition does not truncate the set of possible representations, it
only imposes a restriction on a choice of Euclidean metric on $\R^n$.

\begin{prop}
\label{Euc}
For any linear representation of $\Cliff(r)$ in $\R^n$ there exists a
Euclidean metric on $\R^n$ such that the representation is compatible
with it.
\end{prop}

This fact is well-known but we sketch a proof here.
Consider {\em the Dirac group} consisting of all products of generators $e_1,\dots,e_r$.
Obviously, the Dirac group is finite. Therefore, it has an invariant
Euclidean metric. In particular, all elements $e_i$ act as orthogonal operators
with respect to this metric. A direct computation now shows that all linear
combinations of 1 and the generators act conformally. $\Box$

A representation of a Clifford algebra in a Euclidean space is always assumed
to be compatible with the Euclidean inner product.

Linear representations of the algebras $\C=\Cliff(1)$ and
$\H=\Cliff(2)$ in $\R^n$ are called complex and quaternionic
structures on $\R^n$, respectively. The image of $i$ under a
representation of $\C$ is also called a complex structure. This
may be any linear operator $I$ such that $I^2=-1$. By a complex
(quaternionic) structure on a Euclidean space we always mean a
complex (quaternionic) structure compatible with the Euclidean
metric. For a complex structure $I$ this means that $I$ is
orthogonal and skew-symmetric. For a quaternionic structure, from
the compatibility it follows that all quaternions act conformally,
not only linear combinations of $1,i$ and $j$.

\section{Hopf maps}

We need the following very simple (and well-known) lemma:

\begin{lemma}
\label{cq}
Consider the following fractional linear map from $\A$ to $\A^k$ (where $\A$ is
$\C$ or $\H$):
$$F:x\mapsto(xc+d)^{-1}(xa+b).$$
Here $x,c,d\in\A$ and $a,b\in\A^k$. Assume that the target space is equipped with a
Euclidean metric compatible with the complex or quaternionic structure. Then $F$ takes
all real lines to circles.
\end{lemma}

\proof
We have $F(x)=A+(xc+d)^{-1}B$ where $A=c^{-1}a$ and $B=b-dc^{-1}a$.
Therefore, the image of $\A$ lies on the left line $\{A+uB|\ u\in\A\}$.
The coordinate $u$ on this line is linear and conformal. The map
$x\mapsto u=(xc+d)^{-1}$ is a M\"obius transformation. Hence any real line
goes to a circle. $\Box$

Here is the main construction. Let $V$ be a Euclidean space and $\phi$ a
representation of $\Cliff(V)$ in $\R^n$.
Recall that $\hat V=V\oplus\<1\>\subset\Cliff(r)$.
The following theorem is a partial case of a more general fact,
see e.g. \cite{Yiu}.

\begin{thm}
\label{constr}
Consider the map $F_\phi:(\hat V-0)\times\R^n\to\R^n$ given by
$(\alpha,x)\mapsto\phi^{-1}(\alpha)x$. This map takes all lines to
circles.
\end{thm}

\proof
Consider a line $L$ in $\hat V\times\R^n$. If it passes through 0, then its image is
empty or is a single point. Suppose now that $L$ does not contain the origin.
The projection of $L$ to $\hat V$ generates a subalgebra in $\Cliff(V)$, which is
isomorphic either to $\C$ or to $\H$.

In the first case, $L$ lies in $\C\times\R^n$. The action of $\C$ defines a
complex structure on $\R^n$. Thus $\R^n$ gets identified with $\C^{n/2}$.
Let $t$ be a linear parameter on $L$.
Then the image of $L$ in $\R^n$ is given by $\alpha(t)^{-1}x(t)$ where $x(t)$
is a vector from $\R^n$ and $\alpha(t)$ is a complex number, both depending affinely
on $t$. Now from Lemma \ref{cq} it follows that the image of $L$ is a circle.

The quaternionic case is completely analogous.
$\Box$

{\bf Remark.} The construction given above may be improved. Note that
under the map $F_\phi$ proportional vectors have the same
image. Therefore, there is a map from the projectivization of $\hat V\times\R^n$ to
the sphere $S^n$, which is obtained from $\R^n$ by a one-point compactification.
This map takes projective lines to circles. The projectivization of
$0\times\R^n$ goes to infinity. This is a Hopf map.
The Hopf construction gives the following maps:
$$\RP^3\to S^2,\quad \RP^7\to S^4,\quad \RP^7\to S^6,\quad \RP^{15}\to S^8,\dots$$

\section{Micro-local theorem}

Fix a Euclidean metric on $\R^n$.
We will be assuming that $n\geq 2$.
The inner product of two vectors $x,y\in\R^n$
is denoted by $(x,y)$. Suppose that $\R^n$ is linearly embedded into $\R^{r+n}$.
Let us say that a germ of smooth map $\Phi:(\R^{r+n},0)\to(\R^n,0)$ is a {\em local
projection} if $\Phi$ is the identity on $\R^n$.

\begin{thm}
\label{uni}
Consider a local projection $\Phi:(\R^{r+n},0)\to (\R^n,0)$ that
takes vector lines to circles. Assume that the images of lines not
belonging to $\R^n$ do not degenerate to lines.
Then there is a Euclidean space $V$, a representation $\phi$ of $\Cliff(V)$
in $\R^n$ and an embedding of $\R^{r+n}$ into $V\times\R^n$ identical on $\R^n$
such that the map
$$
F:V\times\R^n\to\R^n,\quad (\alpha,x)\mapsto\phi(1-\alpha)^{-1}x
$$
takes germs of all vector lines to the same circles as $\Phi$ does.
\end{thm}

Theorem \ref{uni} gives a geometric description of representations of Clifford algebras
in terms of circles.

{\bf Remark.} The maps $F$ and $\Phi$ do not need to be the same on $\R^{r+n}$. For instance,
the map $y\mapsto F(\nu(y)y)$, where $\nu$ is an arbitrary smooth function on $\R^{r+n}$,
takes any vector line to the same circle as $F$ does. To make this map satisfy
the conditions of Theorem \ref{uni} we should take a function $\nu$, which is
identically equal to 1 on $\R^n$.

The map $\Phi$ from Theorem \ref{uni} has the form
$$\Phi(\alpha,x)=x+\Gamma(\alpha,x)+\Delta(\alpha)+\dots,$$
where $\Gamma$ is some bilinear map $\R^r\times\R^n\to\R^n$, the form
$\Delta:\R^r\to\R^n$ is quadratic on $\R^r$, and dots denote higher order terms.

\begin{lemma}
\label{div}
The form $\Delta$ vanishes for all $\alpha$.
For any fixed $\alpha$ the quadratic forms $(\Gamma,\Gamma)$ and $(x,\Gamma)$ are
divisible by $(x,x)$.
\end{lemma}

This lemma is an immediate generalization of Proposition 1.3 from \cite{Tim}, and it
can be proved in the same way. We will sketch here another proof assuming for
simplicity that the map $\Phi$ is analytic; see \cite{rounding} for details.

\proof
Extend the Euclidean form $(\cdot,\cdot)$ to $\C^n$ by bilinearity. {\em A complex circle}
in $\C^n$ is by definition the intersection of a sphere
$$\{x\in\C^n|\ (x-a,x-a)=R^2\},\quad a\in\C^n,\ R\in\C$$
with a 2-dimensional complex plane.
A complex circle is an irreducible second degree curve, or a pair of lines
(possibly coincident), or a plane. Obviously, any honest real circle belongs
to a unique complex circle.

The germ $\Phi$ at $0$ extends to the germ of an analytic map from
$(\C^{r+n},0)$ to $(\C^n,0)$. Under restriction to any real vector line,
the function $(\Phi,\Phi)$ and all components of the map $\Phi$
span no more than a 2-dimensional space.
This condition is equivalent to the fact that the image of the line
lies in a circle. On the other hand, it is expressed by a system of
algebraic equations on the Taylor coefficients of the map $\Phi$.
Therefore, the same condition is satisfied under restriction to any
complex line. Thus, the image of the germ of any complex vector line
at 0 lies in a complex circle or in a line.

Denote by $\Nc$ the isotropic cone $\{(x,x)=0\}$ in $\C^n$.
Consider a line $L$ in $\C^{r+n}$ spanned by a vector $(\alpha,x)$ such that
$x\in\Nc$ but $x\ne 0$. The image of $L$ lies in the complex circle $C$
and is tangent to $x$ at 0. We claim that $\Phi(L)$ lies entirely in the
cone $\Nc$. If $C$ is a line or a pair of lines, then $\Phi(L)$ lies in a line,
spanned by $x$ (the algebraic curve $\Phi(L)$ cannot jump from one line to another).
In particular, $\Phi(L)$ lies in $\Nc$.

If $C$ is a plane or an irreducible second degree curve, then it lies entirely
in $\Nc$. Indeed, consider a plane containing the circle $C$.
If this plane does not lie in the cone $\Nc$, then it intersects $\Nc$
in a pair of lines, possibly coincident. The circle $C$ is tangent to
one of these lines at the origin and is asymptotic to both
(i.e. intersects them at infinity), since
the set of infinite points in the projective closure of any complex sphere
is the same as that of $\Nc$. But a second degree curve can not be tangent
to a line at the origin and, at the same time, intersect the same line
at infinity, unless the curve coincides with this line.

Thus for $x\ne 0$ the linear span of the vectors $x$
and $\Phi(\alpha,x)$ lies entirely in $\Nc$. In other words,
$$(x,x)=(x,\Phi(\alpha,x))=(\Phi(\alpha,x),\Phi(\alpha,x))=0$$
for any vector $x\in\Nc$. Therefore, the functions $(x,\Phi(\alpha,x))$
and $(\Phi(\alpha,x),\Phi(\alpha,x))$ are divisible by $(x,x)$ in the ring
of formal power series. Confining ourselves to the consideration of second
order terms only, we obtain the statement of the lemma.
$\Box$

Note that the bilinear map $\Gamma$ contains all information about circles
in the image of $\Phi$. Indeed, to define a circle it is enough to know its
``velocity'' and ``acceleration'' at 0 with respect to some parameter, provided that
the velocity is nonzero. Therefore, if $x\ne 0$, then $\Gamma(\alpha,x)$
uniquely determines the circle that is the image of the line spanned by
$(\alpha,x)$. The subspace $x=0$ goes to the origin under $\Phi$, as it is
seen from the proof of Lemma \ref{div}. This is not difficult to prove
even without the analyticity assumption, with the help of the implicit function theorem.
Thus $\Gamma$ determines the images of all vector lines under $\Phi$.

Theorem \ref{uni} now reduces to a description of the form $\Gamma$.

\section{Connection with the Hurwitz--Radon theorem}

We are going to explore the properties of the bilinear map
$\Gamma:\R^r\times\R^n\to\R^n$ introduced in the previous section.
Recall that $\R^n$ is equipped with a Euclidean metric. Up to now,
$\R^r$ was only a vector space. But we will introduce a Euclidean metric
on it as well. Fix an element $\alpha\in\R^r$. It gives rise to a
linear map $A:\R^n\to\R^n$ given by $A(x)=\Gamma(\alpha,x)$. By Lemma
\ref{div} this map is conformal. This means that the product of $A$
with its conjugate is the multiplication by a real number $q(\alpha)$.
This number is nonzero, since for $x\ne 0$ the image of the line
spanned by $(\alpha,x)$, is an honest circle, in particular,
$\Gamma(\alpha,x)\ne 0$.
It is readily seen that $q$ is a positive definite quadratic form
on $\R^r$. Thus $\R^r$ may be regarded as a Euclidean space.

Now denote by $|\cdot|$ the Euclidean norms in both $\R^r$ and $\R^n$. Then we have
$$|\Gamma(\alpha,x)|=|\alpha|\cdot|x|.$$
Hurwitz \cite{Hu} and Radon \cite{Ra} described all relations of this form.
In terms of Clifford algebras representations their result can be formulated as
follows \cite{ABS}.

\begin{thm}
\label{Clsq}
Let $f:\R^{k+1}\times\R^n\to\R^n$ be a bilinear map such that
$|f(\alpha,x)|=|\alpha|\cdot|x|$ for all $\alpha\in\R^{k+1}$ and $x\in\R^n$.
Then there is a representation $\phi$ of the Clifford algebra $\Cliff(k)$ in
$\R^n$ such that
$$f(\alpha,x)=\phi(\alpha_0+\alpha_1e_1+\cdots+\alpha_ke_k)A_0(x)$$
where $A_0$ is an orthogonal transformation and $\alpha_0,\dots,\alpha_k$ are
coordinates of $\alpha$ in some orthonormal basis.
\end{thm}

In particular, this theorem is applicable to our bilinear map $\Gamma$.
Theorem \ref{Clsq} has a degree of freedom: an arbitrary orthogonal transformation
in the preimage. Our additional restriction on $\Gamma$ kills this degree of freedom.

Some interesting generalizations of Theorem \ref{Clsq}, and the general Hurwitz
problem are discussed in \cite{Sh}.

\section{Complex multiplications and Clifford algebras}

In this section, we establish a connection between the bilinear map $\Gamma$ and
a representation of some Clifford algebra. This connection is in the spirit
of Theorem \ref{Clsq}, and we use similar arguments.

Fix a Euclidean metric $(\cdot,\cdot)$ on $\R^n$.
A linear operator $A:\R^n\to\R^n$ is called a {\em complex multiplication} if
either $A$ is the multiplication by a real number, or there exists a
complex structure $I$ on $\R^n$ compatible with $(\cdot,\cdot)$ such that $A$ is
the multiplication by a complex number with respect to $I$.

Here is a characterization of all complex multiplications in real terms.

\begin{prop}
\label{Cm}
A linear operator $A:\R^n\to\R^n$ is a complex multiplication if and only if
both quadratic forms $(x,Ax)$ and $(Ax,Ax)$ are divisible by $(x,x)$.
\end{prop}

\proof
Any complex structure compatible with the Euclidean metric is represented
by an operator $I$, which is orthogonal and skew-symmetric at the same time.
In other words, the equalities $(Ix,Ix)=(x,x)$ and $(Ix,x)=0$ hold for
each $x\in\R^n$. The second equality follows from the first equality
and the relation $I^2=-1$.
If $A=a+bI$ for some complex structure $I$, then, as it is readily seen,
$(x,Ax)=a(x,x)$ and $(Ax,Ax)=(a^2+b^2)(x,x)$.

Assume now that $(x,Ax)=p(x,x)$ and $(Ax,Ax)=q(x,x)$ for some real
numbers $p$ and $q$.
From the Cauchy-Schwartz inequality $(Ax,Ax)(x,x)\geq (x,Ax)^2$
it follows that $q\ge p^2$. If $q=p^2$,
then $A$ is everywhere proportional to $x$.
Since $(x,Ax)=p(x,x)$, the coefficient of proportionality equals to $p$.
Thus $A$ is the multiplication by $p$. If $q>p^2$, then $A=p+\sqrt{q^2-p^2}I$
for some linear operator $I$.
From the restrictions on $A$ it follows that $I$
is both orthogonal and skew-symmetric.
The orthogonality means that $I^*I=1$ where $I^*$ is the conjugate
operator to $I$. The skew-symmetry says that $I^*=-I$.
Hence $I^2=-1$, i.e. $I$ is a complex structure, and $A$ is a complex
multiplication. $\Box$

Let $\Gamma:\R^r\times\R^n\to\R^n$ be the bilinear map introduced in Section 4.
Identify $\R^r$ with a space of complex multiplications: a vector $\alpha\in\R^r$
gets identified with the operator $x\mapsto\Gamma(\alpha,x)$, which is a complex
multiplication by Lemma \ref{div} and Proposition \ref{Cm}. Denote by $\hat V$ the
linear span of $\R^r$ and the identity transformation.
There is a Euclidean form on $\hat V$. Namely, for each complex multiplication
$A$ put $q(A)=(A(x),A(x))/(x,x)$, which is independent of $x$. Clearly, $q$ is
a positive definite quadratic form.

\begin{prop}
\label{Cl}
Let $V$ be the orthogonal complement to $1$ in $\hat V$ with respect to $q$.
Then there exists a representation $\phi$ of $\Cliff(V)$ in $\R^n$ such that for
any $\alpha\in\R^r$ and $x\in\R^n$ we have $\Gamma(\alpha,x)=\phi(\alpha)(x)$
under the natural embedding of $\hat V$ into $\Cliff(V)$.
\end{prop}

\proof
Denote by $q(\cdot,\cdot)$ the polarization of the quadratic form $q(\cdot)$.
By definition, for any $A\in\hat V$ and $x\in\R^n$ we have $(Ax,Ax)=q(A)(x,x)$.
Therefore, any 2 elements $A,B\in\hat V$ satisfy the relation
$(Ax,By)+(Bx,Ay)=2q(A,B)(x,y)$
identically with respect to $x,y\in\R^n$. In other words, $A^*B+B^*A=2q(A,B)$.

Choose an orthonormal basis in $\hat V$ containing 1. Denote the other elements
of this basis by $I_1,\dots,I_r$. We have
$$I_j+I_j^*=2q(1,I_j)=0,\quad I_j^2=-I_j^*I_j=-q(I_j)=-1,$$
$$I_jI_k+I_kI_j=-(I_j^*I_k+I_k^*I_j)=-q(I_j,I_k)=0.$$
Thus the operators $I_1,\dots,I_r$ give rise to a representation of $\Cliff(V)$.
$\Box$

{\sc Proof of Theorem \ref{uni}.}
By Proposition \ref{Cl}, the bilinear map $\Gamma$ sends $(\alpha,x)$ to
$\phi(\alpha)(x)$
for some embedding of $\R^r$ into $\hat V$ and some representation $\phi$ of $\Cliff(V)$
in $\R^n$. Thus the 2-jet of $\Phi$ at 0 is the same as that of
$F:(\alpha,x)\mapsto\phi(1-\alpha)^{-1}x$. But the circles in the image are
determined by the 2-jet only. Hence the images of all vector lines under the maps
$\Phi$ and $F$ are the same. $\Box$

Keywords: line, circle, Clifford algebra, Hopf map.

\end{document}